\documentclass[11pt,reqno,oneside]{amsart}

\usepackage{amssymb}
\usepackage{graphics}
\usepackage[usenames]{color}
\usepackage[hidelinks]{hyperref}

\renewcommand{\theequation}{\arabic{section}.\arabic{equation}}

\newtheorem{theorem}{Theorem}[section]
\newtheorem{lemma}{Lemma}[section]

\newtheorem{rem}{Remark}[section]
\newtheorem{defi}{Definition}[section]

\def\R{{\rm I}\!{\rm R}}

\newcommand{\bit}{\begin{itemize}}
\newcommand{\eit}{\end{itemize}}
\newcommand{\bq}{\begin{equation}}
\newcommand{\eq}{\end{equation}}

\newcommand{\ds}{\displaystyle}
\newcommand{\hphi}{\widehat{\phi}}
\newcommand{\hx}{\hat x}

\newcommand{\at}{\tilde a}
\newcommand{\bt}{\tilde b}
\newcommand{\hk}{\widehat K}
\newcommand{\hT}{\widehat T}

\newcommand{\oK}{\overline K}
\newcommand{\oQ}{\overline Q}

\newcommand{\ou}{\overline u}
\newcommand{\ov}{\overline v}
\newcommand{\ox}{\overline x}

\newcommand{\q}{\mathcal Q}

\makeatletter
\g@addto@macro{\endabstract}{\@setabstract}
\newcommand{\authorfootnotes}{\renewcommand\thefootnote{\@fnsymbol\c@footnote}}%
\makeatother

\author[G. Monz\'on]{Gabriel Monz\'on}
\address{Instituto de Ciencias\\ Universidad Nacional de General Sarmiento\\
	J. M. Guti\'errez 1150 \\ 
	(1613) Los Polvorines\\
	Buenos Aires\\ Argentina.}
\email{gmonzon@campus.ungs.edu.ar}
\thanks{The author thanks G. Acosta for his invaluable help and support, and the referees for several suggestions and contributions.}

\title[Anisotropic interpolation error estimate for quadrilateral elements]{Anisotropic interpolation error estimate for arbitrary quadrilateral isoparametric elements}

\begin{document}
\maketitle
\begin{abstract}
	The aim of this paper is to show that, for any $p \in [1,\infty)$, the $W^{1,p}$-anisotropic interpolation error estimate holds on quadrilateral isoparametric elements verifying the \emph{maximum angle condition} ($MAC$) and the \emph{property of comparable lengths for opposite sides} ($clos$), i.e., on all those quadrilaterals with interior angles uniformly bounded away from $\pi$ and with both pairs of opposite sides having comparable lengths. 
		
	For rectangular elements our interpolation error estimate agrees with the usual one whereas for perturbations of rectangles (the most general quadrilateral elements previously considered as far as we know) our result has some advantages with respect to the pre-existing ones: the interpolation error estimate that we proved is written by using two neighboring sides of the element instead of the sides of the unknown perturbed rectangle and, on the other hand, conditions $MAC$ and $clos$ are much simpler requirements to verify and with a clearer geometrical sense than those involved in the definition of perturbations of a rectangle.
\end{abstract}

\section{Introduction}
\label{intro}
\setcounter{equation}{0}

Let $K\subset \R^2$ be an arbitrary convex quadrilateral and $\hk :=[0,1]^2$ the reference unitary square. For $1\le i\le 4$, we denote with $V_i$ (resp. $\widehat{V}_i$) the vertices of $K$ (resp. $\hk$) enumerated in counterclockwise order. In this context, for any $K$, $V_1$ is arbitrarily chosen with the exception of the reference element for which $\hat{V}_1$ is assumed to be placed at the origin. Associated to  each $\widehat{V}_i$  we consider the respective bilinear basis function $\hphi_i$ and define the  mapping $F_K: \hk \to K$,
$F_K(\hx)=\sum_{i=1}^4 V_i \hphi_i(\hx)$. The $\q_1$-isoparametric Lagrange 
interpolation $Q$ is then defined on $K$ by \cite{CR}
$$Qu(x)=\widehat{Q} \hat{u} (\hx),$$

where $x=F_K(\hx)$ and $\widehat{Q}$ is the $\q_1-$ Lagrange interpolation of 
$\hat{u}=u \circ F_K$ on $\hk$.

The classical $H^1$-error estimate \cite{BS,CR} for $Q$ can be written as 
\bq
\label{eq:claserror}
|u - Q u|_{1,K} \le  C h |u|_{2,K}\ , 
\eq
where $h$ is the diameter of $K$ and $C$ is a positive constant.  

The convexity of $K$ is not a sufficient condition in order to guarantee \eqref{eq:claserror}. That is, there are certain families of convex quadrilaterals
for which it is not possible to take a fixed constant $C$ in \eqref{eq:claserror} (see, for intance, Counterexamples 6.1 and 6.2 in \cite{AM}). As a consequence, some extra geometric assumptions on $K$ have been introduced along the past decades. In \cite{CRpaper} it was considered the existence of a constant $\beta$ such that
\bq
\label{cond:CR1}
0 < \beta \le \frac{s}{h} \sqrt{1-\gamma}
\eq
where $s$ is the shortest side of $K$ and $\gamma = \max \{ |\cos(\theta)| : 
\theta \mbox{ is an inner angle of } K \}$. Condition \eqref{cond:CR1} can be 
rewritten \cite{CR} as follows: there exists positive constants $\mu_1, \mu_2$ such that
\bq
\label{cond:CR2}
\frac{h}{s} \le \mu_1 \quad \mbox{ and } \quad 
\gamma \le \mu_2 <1.
\eq
Roughly speaking, the inequality on the left in \eqref{cond:CR2} says that lengths of the four sides of $K$ should be nearly equal, while the inequality on the right provides that each internal angle of $K$ must be bounded away from zero and $\pi$.

In \cite{J}, the author considers the usual {\emph{regularity condition}}: the 
existence of a constant $\sigma$ such that
\bq
\label{eq:rc}
h/\rho \leq \sigma
\eq
where $\rho$ denotes the diameter of the maximum circle contained in $K$. Elements satisfying (\ref{eq:rc}), called  {\it regular elements} in the sequel, can not be too narrow. In spite of its simplicity and the fact that regular quadrilaterals can degenerate into triangles, \eqref{eq:rc} is an undesirable requirement in some particular applications (see \cite{ApD,Ap,Apel,D} and references there in). 

Certain class of quadrilaterals, for which condition \eqref{eq:rc} does not hold, and yet \eqref{eq:claserror} remains valid, was described in \cite{ZV1,ZV2}. Essentially, these quadrilaterals must have their two longest sides as opposite and parallels (or almost parallels) to each other. More recently, was proved \cite{AM1, AM} that only the term on the right in \eqref{cond:CR2} is a simple sufficient condition for \eqref{eq:claserror}. As we pointed out, this condition requires  all internal angles of $K$ to be bounded away from zero and $\pi$. This condition allows to consider a large family of arbitrarily narrow (often called {\emph{anisotropic}}) elements. 

More formally, we introduce some notation in the following definition.
\begin{defi}
	\label{defidac}
	We say that a quadrilateral $K$ satisfies the {\slshape{maximum angle condition}} with constant $\psi_M$, or shortly $MAC(\psi_M)$, if all inner 
	angles $\theta$ of $K$ verify $$\theta \leq \psi_M<\pi.$$ Similarly, we say 
	that a quadrilateral $K$ satisfies the {\slshape{minimum angle condition}} with constant $\psi_m$, or shortly $mac(\psi_m)$, if all inner angles $\theta$ of 	$K$ verify $$0<\psi_m \leq \theta.$$
	
	Finally, we say that a quadrilateral $K$ satisfies the {\slshape{double angle condition}} with constants $\psi_m, \psi_M$, or shortly $DAC(\psi_m,\psi_M)$, 
	if $K$ simultaneously verifies $MAC(\psi_M)$ and $mac(\psi_m)$.
\end{defi}

In this paper we show that if $K$ is an arbitrary convex quadrilateral satisfying the $MAC$, under the extra assumption that any pair of its opposite sides have comparable lengths, then the estimate \eqref{eq:claserror} can be improved by the following {\it anisotropic interpolation error estimate} 
\bq
\label{ep:aeeK}
|u-Qu|_{1,K} \le C \left[ 
|l_1| \left\| \partial_{l_1} \nabla u \right\|_{0,K} + 
|l_2| \left\| \partial_{l_2} \nabla u \right\|_{0,K} 
\right]
\eq
where $l_1$ and $l_2$ are two neighboring sides of $K$ with the property that the parallelogram determined by them contains $K$ entirely (as usual, $\partial_{l_i}$ denotes the directional derivative along $l_i$). 

For future references we formalize the extra assumption on $K$ described before as follows

\begin{defi}
	\label{def:clos}
	Let $K$ be an arbitrary convex quadrilateral. We say that $K$ verifies the {\rm property of comparable lengths for opposite sides} with constants $C_*$ and $C^*$, or shortly $clos(C_*,C^*)$, if for any two opposite sides $l_1$ and $l_2$ of $K$ holds true that
	$$C_* \le \frac{|l_1|}{|l_2|} \le C^*.$$
\end{defi}

A striking fact is that the requirements $MAC$ and $clos$ imply the $mac$ (see Lemma \ref{lem:red} below); and hence, conditions $[MAC,clos]$ are equivalent to $[DAC,clos]$. Keeping this in mind, the extra assumption $clos$ allows to improve the classical interpolation error estimate noticeably. Moreover, under the same hypothesis ($MAC$ and $clos$, eq. $DAC$ and $clos$), a similar estimate to \eqref{ep:aeeK} can be obtained for any $p \in [1,\infty)$. This is, for any $p \in [1,\infty)$, there exists a positive constant $C$ (depending only on those constants involved in the $MAC$ and the $clos$) such that
\bq
\label{ep:aeeKp}
|u-Qu|_{1,p,K} \le C \left[ 
|l_1| \left\| \partial_{l_1} \nabla u \right\|_{0,p,K} + 
|l_2| \left\| \partial_{l_2} \nabla u \right\|_{0,p,K} 
\right].
\eq

For any fixed class of elements containing anisotropic quadrilaterals, the 
estimate \eqref{ep:aeeK} is better than \eqref{eq:claserror}, since  the former allows to keep the error at a fixed level by selecting appropriately the size of the mesh in each independent direction rather than enforcing an uniform reduction of both of them.

Anisotropic error estimates for rectangles or parallelograms (affine images of the unitary square) can be found in  \cite{ApD,Ap}; however, for general quadrilateral elements, the bibliography known is scarce. To our best 
knowledge the main results about this topic is stated in \cite{Apel}. 
Basically, quadrilaterals considered in that work are perturbations of a rectangle in the following sense: (after rigid movements we can assume that the element has the origin as one of its vertices) there exists positive constants $a$ and $b$, with $b \le a$, such that the map $F_K$ between the unitary square $\hk$ 
and the element $K$ can be written as
\bq
\label{eq:mapapel}
\ds F_K(\hx) = (a\hx_1, b\hx_2) + \sum_{i=1}^4 a^{(i)} \hphi_i(\hx) 
\eq
where $\hphi_i$ denotes, as before, the basis nodal function associated to 
$\widehat{V}_i$ on $\hk$, and for the distortive vectors $a^{(i)}=(a^{(i)}_1, 
a^{(i)}_2)$, there exists constants $a_0,a_1, a_2$ such that
\bq
\label{eq:vd}
|a_i^{(j)}| \leq a_i b, \quad 0 \leq a_i \lesssim 1, \quad i=1,2,\ j = 1, 2, 3, 4,
\eq 
and
\bq
\label{eq:jac}
\ds \frac{1}{2}-\frac{b}{a}a_1-a_2 \geq a_0 > 0.
\eq

For any rectangle perturbation $K$, it was proved \cite[Theorem 2.8]{Apel} the following anisotropic error estimate
\bq
\label{eq:aeeapel}
|u - Qu|_{1,K} \le C \left[ 
a \left\| \frac{\partial}{\partial x_1} \nabla u \right\|_{0,K} + 
b \left\| \frac{\partial }{\partial x_2} \nabla u \right\|_{0,K} \right].
\eq

We believe that  \eqref{ep:aeeK} has some advantages over
\eqref{eq:aeeapel}. Indeed,  the former is stated in terms of the sides of the quadrilateral $K$ itself instead of  the sides of
the rectangle from which $K$ is obtained after a perturbation. Actually, although this rectangle  can be easily found from 
the geometry of $K$,  it is  not explicitly given. On the other hand, conditions \eqref{eq:vd} and \eqref{eq:jac} have not a clear geometrical interpretation  
whereas conditions $MAC$ and $clos$ do have it, which can be useful in practical tests.

\medskip
The outline is as follows: In Section \ref{impli} we introduce a particular class of elements called {\it the reference setting} and we give a characterization for all those quadrilaterals satisfying the $MAC$ and the $clos$ in terms of these reference elements. Section \ref{sub:red} is devoted to prove that our study can be reduced to the reference setting; that is, to prove that if the anisotropic interpolation error estimate holds on the reference elements, then the anisotropic interpolation error estimate holds on arbitrary quadrilaterals. Sections \ref{errtreat} and \ref{ontriangles} are dedicated to show the error treatment and to present appropiate bounds for each term appearing in such decomposition. Finally, in Section \ref{main} we present our main results.

\section{Implications of the $DAC$ and the $clos$}
\label{impli}
\renewcommand{\theequation}{\arabic{section}.\arabic{subsection}.\arabic{equation}}
\setcounter{equation}{0}

\subsection{Basic notation and some useful facts}
\setcounter{equation}{0}

Given $a,b,\at,\bt>0$, $K(a,b,\at,\bt)$ will represent a convex quadrilateral with vertices $V_1=(0,0)$, $V_2=(a,0)$, $V_3=(\at,\bt)$ and $V_4=(0,b)$. We usually refer to this kind of quadrilateral as a {\it reference element}.

\begin{figure}[h]
	\centering
	\resizebox{15cm}{7cm}{\includegraphics{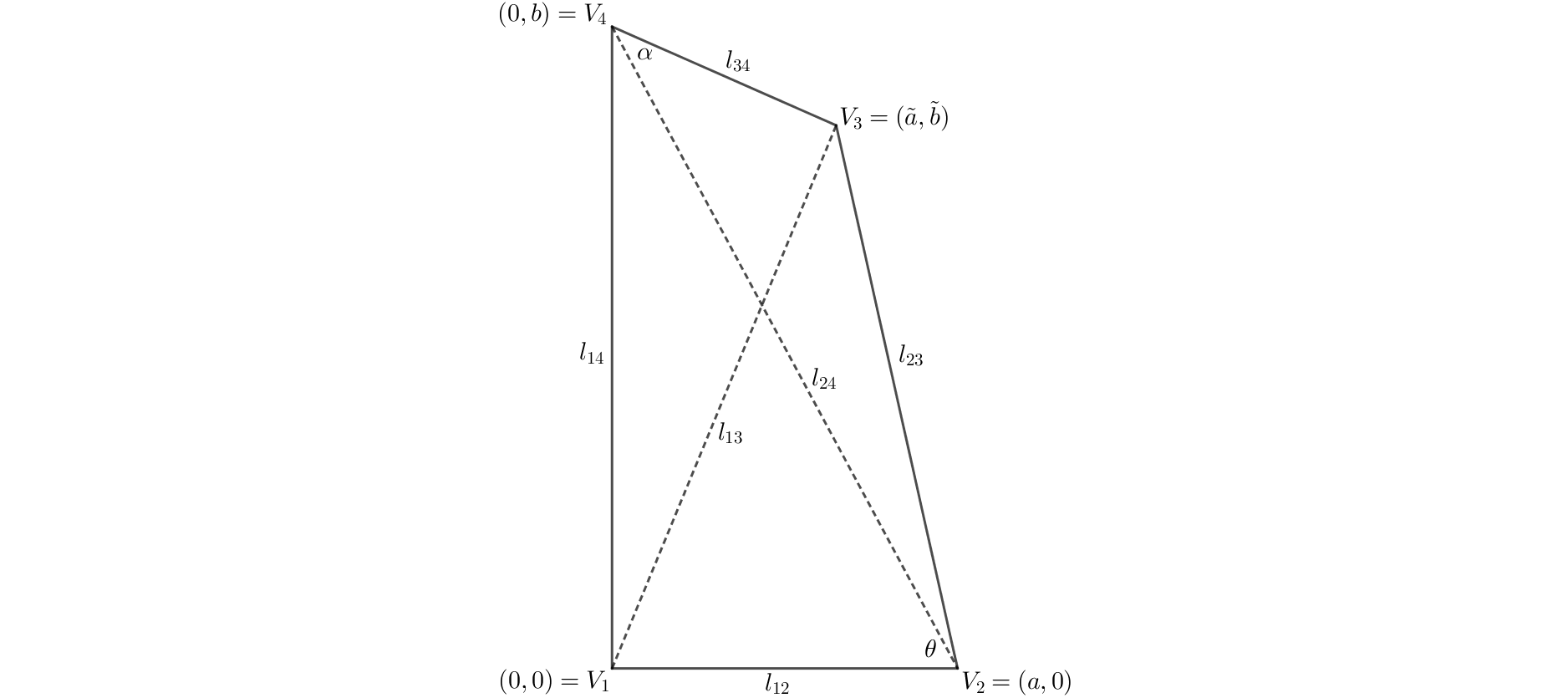}}	
	\caption{Basic notation adopted on a reference element $K(a,b,\at,\bt)$.}
	\label{fig:not}
\end{figure}

The segment connecting $V_i$ to $V_j$ will be denoted by $l_{ij}$ (see Figure \ref{fig:not}).

We will be particularly interested on those reference elements $K(a,b,\at,\bt)$ that verify conditions $(D1)$ and $(D2)$ given by \eqref{cond:D1} and \eqref{cond:D2}, respectively. Roughly speaking, the condition $(D1)$ says that $V_3=(\at,\bt)$ is an interior point of the rectangle $R_{ab}=K(a,b,a,b)$, indeed
\bq
\label{cond:D1}
\mbox{$(D1)$} \qquad \frac{\at}{a},\ \frac{\bt}{b} \le 1.
\eq 

On the other hand, the condition $(D2)$ ensures that the angle $\alpha$ of $\Delta(V_2V_3V_4)$ placed at $V_4$ (see Figure \ref{fig:not}) is bounded away from zero and $\pi$. This is, there exists a positive constant $C$ such that
\bq
\label{cond:D2}
\mbox{$(D2)$} \qquad \frac{1}{\sin(\alpha)} \le C.
\eq

Borrowing the convention assumed in \cite{AM}, we say that quadrilaterals 
$K_1$ and $K_2$ are {\it equivalents} if there exists an invertible affine mapping $L$ such that $L(K_1)=K_2$ and the matrix 
$B$ associated to $L$ satisfies $\|B\|,\|B^{-1}\| \le C$ for some constant $C$. In this case, we write $K_1 \equiv_C K_2$.

Notice that, for a such matrix $B$, it holds that $\kappa(B) \le C^2$ where $\kappa(B)$ denotes the condition number of $B$. 

\medskip
We conclude this brief section by showing that both conditions $DAC$ and $clos$ are inherited properties between equivalent elements.

\begin{lemma}
\label{lem:closueq}	
Let $K_1$ and $K_2$ be two convex quadrilaterals such that $K_1 \equiv_C K_2$. If $K_1$ verifies the $clos(C_*,C^*)$, then $K_2$ verifies the $clos(C^{-2}C_*,C^2C^*)$. In particular, $K_1$ satisfies the $clos$ if and only if $K_2$ satisfies the $clos$. 
\end{lemma}

\begin{proof}
Since $K_1 \equiv_C K_2$, there exists an invertible linear mapping $L: K_1 \to K_2$ given by $L(x)=Bx+P$ such that $\|B\|,\|B^{-1}\| \le C$.

Let $e$ be an arbitrary side of $K_1$. Without loss of generality, we can assume that $e$ is the segment joining $V_1$ and $V_2$. Then
$$|L(e)| = \left\| L(V_2) - L(V_1) \right\| = 
\left\| B(V_2-V_1) \right\| \le 
\left\| B \right\| \left\| V_2-V_1 \right\| = 
\left\| B \right\| |e| \le
C|e|$$
and, on the other hand, 
$$|e| = \left\| L^{-1}(L(V_2)) - L^{-1}(L(V_1)) \right\| = 
\left\| B^{-1}(L(V_2)-L(V_1)) \right\| \le 
\left\| B^{-1} \right\| |L(e)| \le
C|L(e)|.$$
Therefore, 
\bq
\label{eq:osc}
C^{-1}|L(e)| \le |e| \le C |L(e)|.
\eq

The proof concludes immediately taking into account that $L$ applies opposite sides of $K_1$ into opposite sides of $K_2$ and using \eqref{eq:osc}. \qed

\end{proof}

\medskip
In \cite{AD1} it was proved the following elementary result about the behavior of angles under an affine mapping.

\begin{lemma}
	\label{lema:angles}
	Let $L$ be an affine transformation associated with a matrix $B$. Given two 
	vectors 
	$v_1$ and $v_2$, let $\alpha_1$ and $\alpha_2$ be the angles between them and 
	between $L(v_1)$ and $L(v_2)$, respectively. Then
	\bq
	\label{eq:angles}
	\frac{2}{\kappa(B)\pi} \alpha_1 \le \alpha_2 \le \pi \left( 
	1-\frac{2}{\kappa(B)\pi} \right) + \alpha_1 \frac{2}{\kappa(B)\pi}.
	\eq
\end{lemma}

\begin{proof}
	We refer to \cite[Lemma 5.6]{AD1} for the proof. \qed
\end{proof}

\medskip
If $\alpha_1$ is bounded away from zero and $\pi$, from \eqref{eq:angles}, it follows that $\alpha_2$ does not approach zero nor $\pi$. This simple fact allow us to prove the following result, however we omit the details of its proof.

\begin{lemma}
	\label{lem:DACueq}
	If $K_1$ and $K_2$ are two convex quadrilaterals such that $K_1 \equiv_C K_2$, then $K_1$ verifies the $DAC(\psi_m,\psi_M)$ if and only if $K_2$ verifies the $DAC(\overline{\psi}_m,\overline{\psi}_M)$. 
\end{lemma}

In Lemma \ref{lem:DACueq}, the constants $\overline{\psi}_m$ and $\overline{\psi}_M$ depend on $\psi_m$, $\psi_M$ and $C$; and, reciprocaly, the constants $\psi_m$ and $\psi_M$ depend on $\overline{\psi}_m$, $\overline{\psi}_M$ and $C$. The explicit dependency between all these constants is not necessary for our purposes, so we omit it, but it can be easily deduced from \eqref{eq:angles}.

\subsection{Implications of the $DAC$}
\setcounter{equation}{0}

Theorem 3.3 in \cite{AM} says that any convex quadrilateral satisfying the $DAC$ is equivalent to a reference element $K(a,b,\at,\bt)$ obeying $[D1,D2]$ (given by \eqref{cond:D1} and \eqref{cond:D2}, respectively). We rewrite this result in Theorem \ref{teo:carDAC} below and, for the sake of completeness, we reproduce its proof which will be also useful to clarify some future references and notation. 

\begin{figure}[h]
	\centering
	\resizebox{15.5cm}{7cm}{\includegraphics{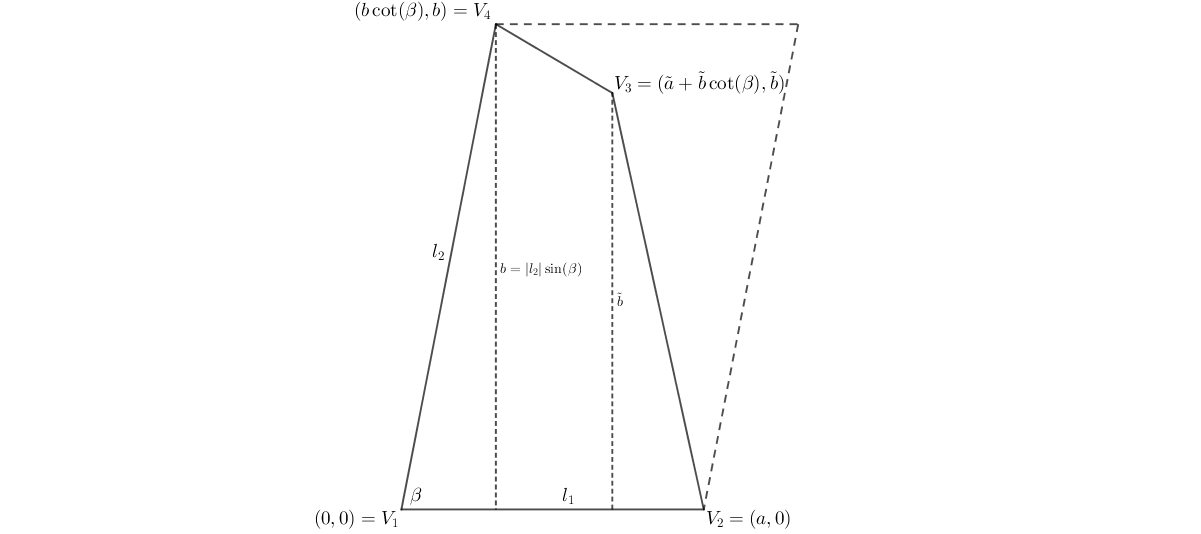}}	
	\caption{Notation adopted on a general quadrilateral $K$.}
	\label{fig:notonK}
\end{figure}

\begin{theorem}
\label{teo:carDAC}
Let $K$ be a general convex quadrilateral. Then $K$ satisfies the $DAC(\psi_m,\psi_M)$ if and only if $K$ is equivalent to an element $K(a,b,\at,\bt)$ verifying $[D1,D2]$, where constant $C$ involved in $(D2)$ only depends on $\psi_m$ and $\psi_M$.
\end{theorem}

\begin{proof}
	It is always possible to select, maybe not in a unique way, two neighboring sides $l_1,l_2$ of $K$ such that the parallelogram defined by them contains $K$ entirely. Making a choice if necessary, we choose $l_1$ and $l_2$ with the properties just described. Call $V_1$ the common vertex of $l_1$ and $l_2$, and let $\beta$ be the angle placed at $V_1$. After a rigid movement, we may assume that $V_1=(0,0)$ and that $l_1=V_1V_2$ lies along the $x$ axis (with nonnegative coordinates), i.e., the vertex $V_2$ is placed at $(a,0)$ for some $a>0$. Moreover, we can also assume that $l_2$ belongs to the upper half plane. Since $l_2=V_1V_4$ we get that $V_4=(b \cot(\beta),b)$ with $b = |l_2| \sin(\beta)$ (see Figure \ref{fig:notonK}). We affirm that the linear mapping $L$ associated to the matrix $B=\left( \begin{array}{cc} 1&\cot(\beta)\\ 0&1\end{array} \right)$ performs the desired transformation. Indeed, since  $\|B\|,\|B^{-1}\|<\frac{\sqrt{2}}{\sin(\beta)} \le C$ with $C=C(\psi_m,\psi_M)$ (thanks to the assumption $DAC$ on $K$) and choosing $\at, \bt$ in such a way that $L(\at,\bt)=V_3$ ($V_3$ denotes the remaining vertex of $K$), it follows that $L$ belongs to the class of affine transformations involved in the definition of {\it equivalent} quadrilaterals. Clearly, $(D1)$ holds. On the other hand, thanks to Lemma \ref{lema:angles} and the fact that the interior angle of $K$ placed at $V_3$ is away from $0$ and $\pi$ by the assumption $DAC$, the interior angle of $K(a,b,\at,\bt)$ placed at $(\at,\bt)$ is also bounded away from zero and $\pi$ meaning that at least one of the remaining angles of the triangle of vertices $(0,b),(\at,\bt),(a,0)$ does not approach zero or $\pi$. Performing a rigid movement if necessary, we may assume that this is the one at $(0,b)$ and hence $(D2)$ follows.  
	Reciprocally, assume that $K(a,b,\at,\bt)$ verifies $[D1,D2]$ and it is equivalent to $K$. Notice that maximal and minimal angle of $K(a,b,\at,\bt)$ are away from $\pi$ and $0$, respectively (in terms of those constants involved in $[D1,D2]$). Indeed, since at $V_1$ we have a right angle we only need to check the remaining vertices. The angle placed at $V_4$ is bounded above by $\pi/2$ due to $(D1)$ and below by $\alpha$. Let us focus now on the angle at vertex $V_3$. It should be bounded below by $\pi/2$ due to $(D1)$. On the other hand, it can not approach to $\pi$ due to $(D2)$. Finally, the angle at $V_2$ is greater than $\alpha$ and also bounded above by $\pi/2$. The proof concludes by using that $K$ is equivalent to $K(a,b,\at,\bt)$ and Lemma \ref{lem:DACueq}. \qed
\end{proof}

\begin{figure}[h]
	\centering
	\resizebox{16cm}{7cm}{\includegraphics{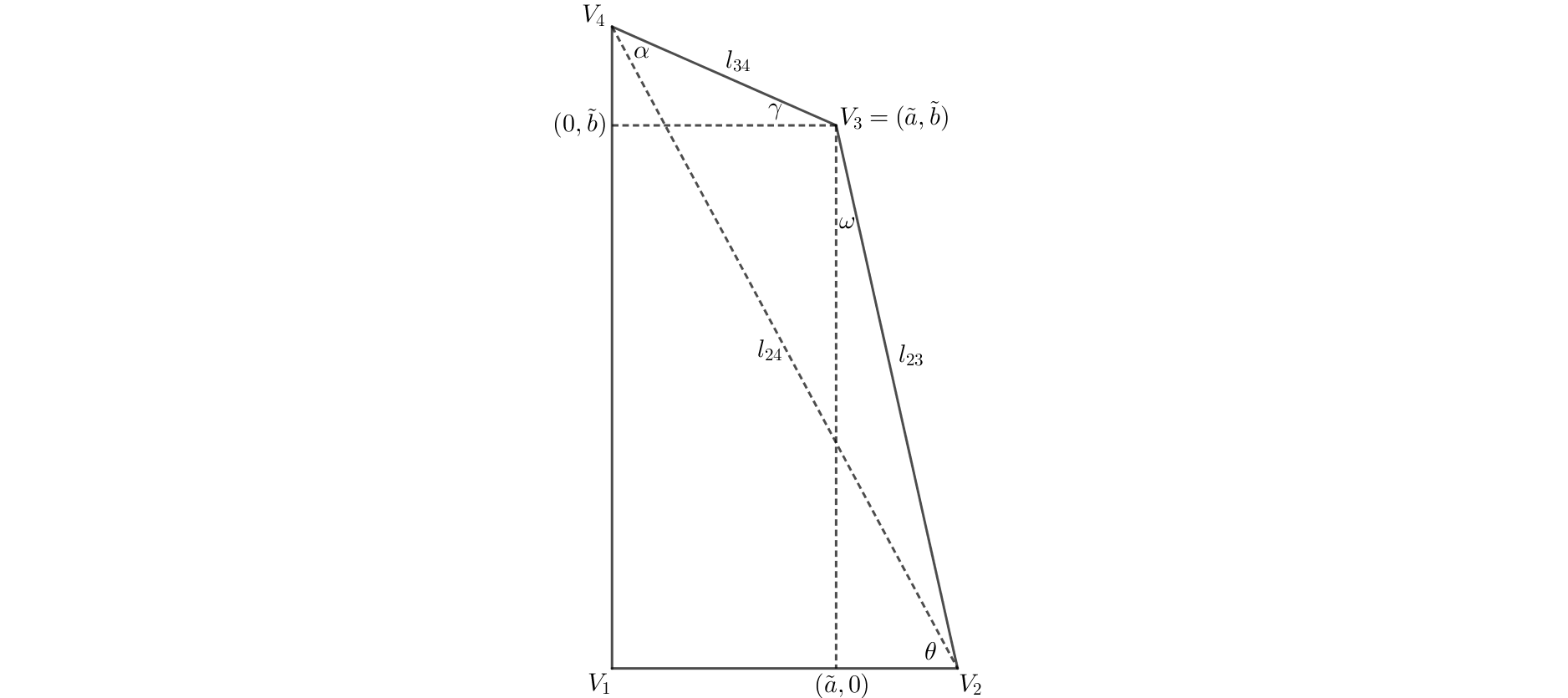}}	
	\caption{Notation on a reference element $K(a,b,\at,\bt)$ satisfying conditions $[D1,D2]$.}
	\label{fig:rem}
\end{figure}
\medskip
We remark some elementary but useful facts on a convex quadrilateral $K=K(a,b,\at,\bt)$ obeying conditions $[D1,D2]$.

\begin{rem}
	The diameter $h$ of $K$ is attained on the diagonal that connects $V_2=(a,0)$ with $V_4=(0,b)$ $($see {\rm Figure \ref{fig:rem})}, i.e.
	\bq
	\label{heqv2v4}
	h=|l_{24}|.
	\eq
\end{rem}  

\begin{rem}
	There exists a constant $C$ such that
	\bq
	\label{aleb}
	a \le C b.
	\eq
	Actually, $C$ can be taken as the same constant appearing in $(D2)$. Indeed, let $\theta$ be the interior angle of $\Delta(V_1V_2V_4)$ placed at $V_2$. It is clear that $\alpha \le \theta < \pi/2$ $($see {\rm Figure \ref{fig:rem})}, then $($thanks to $(D2))$
	$$\frac{a}{b}=\cot(\theta) \le \cot(\alpha) \le C.$$
	
	In particular, notice that the case $a<<b$ is allowed.
\end{rem} 

\begin{lemma}
	\label{lem:D1eqclos}
	Let $K(a,b,\at,\bt)$ be a convex quadrilateral satisfying $[D1,D2]$. Then, $K(a,b,\at,\bt)$ obeys the $clos$ if and only if there exists a positive constant $C_0 < 1$ such that
	\bq
	\label{cond:D3}
	\mbox{$(D3)$} \qquad C_0 \le \frac{\at}{a},\ \frac{\bt}{b}.
	\eq 
\end{lemma}

\begin{proof}
	Let $\gamma$ be the angle between $l_{34}$ and the segment that connects $V_3$ with $(0,\bt)$, and let $\omega$ be the angle between $l_{23}$ and the segment that connects $V_3$ and $(\at,0)$, respectively (see Figure \ref{fig:rem}). Notice that $K(a,b,\at,\bt)$ satisfies the $DAC(\psi_m,\psi_M)$ for some constants $\psi_m,\psi_M$ (thanks to Theorem \ref{teo:carDAC}); therefore, $\gamma+\omega \le \psi_M-\pi/2$ since the interior angle $\beta_3=\gamma+\omega+\pi/2$ of $K$ placed at $V_3$ is bounded above by $\psi_M$. Then, it follows that $\gamma, \omega \le \psi_M - \frac{\pi}{2} < \frac{\pi}{2}$ and hence
	\bq
	\label{tangamma}
	 \frac{b-\bt}{\at} = \tan(\gamma) \le  \tan(\psi_M-\pi/2)
	\eq
	\bq
	\label{tanomega}
	\frac{a-\at}{\bt} = \tan(\omega) \le  \tan(\psi_M-\pi/2).
	\eq
	
	Since $\ds \frac{|l_{34}|^2}{|l_{12}|^2} = \frac{\at^2+(b-\bt)^2}{a^2}$, from \eqref{tangamma} and thanks to the non-negativity of $b-\bt$, we obtain
	\bq
	\label{horiz}
	\frac{\at^2}{a^2} \le \frac{|l_{34}|^2}{|l_{12}|^2} \le (1+\tan^2(\psi_M-\pi/2)) \frac{\at^2}{a^2}.
	\eq
	
	Similarly, from the fact $\ds \frac{|l_{23}|^2}{|l_{14}|^2} = \frac{(a-\at)^2+\bt^2}{b^2}$, \eqref{tanomega} and  the non-negativity of $a-\at$, we obtain
	\bq
	\label{vert}
	\frac{\bt^2}{b^2} \le \frac{|l_{23}|^2}{|l_{14}|^2} \le (1+\tan^2(\psi_M-\pi/2)) \frac{\bt^2}{b^2}.
	\eq
	
	Finally, taking into account that $\at/a, \bt/b \le 1$ (due to $(D1)$), Lemma follows easily from \eqref{horiz} and \eqref{vert}. \qed
\end{proof}

\subsection{Implications of the $DAC$ and the $clos$} 
\setcounter{equation}{0}

We are interested on those elements that satisfy the $DAC$ and the $clos$. The following theorem provides a characterization of such quadrilaterals in terms of the reference elements that will be useful to our purposes.

\begin{theorem}
	\label{teo:DAC+clos}
	Let $K$ be a general convex quadrilateral. Then $K$ satisfies the $DAC$ and the $clos$ if and only if $K$ is equivalent to an element $K(a,b,\at,\bt)$ verifying conditions $[D1,D2,D3]$ $($given by \eqref{cond:D1}, \eqref{cond:D2} and \eqref{cond:D3}, respectively$)$.
\end{theorem}

\begin{proof}
	Assume that $K$ satisfies the $DAC$ and the $clos$. From Theorem \ref{teo:carDAC}, $K$ is equivalent to an element $K(a,b,\at,\bt)$ satisfying conditions $[D1,D2]$. Since the $clos$ is a property that is preserved between equivalent elements (Lemma \ref{lem:closueq}), it follows that $K(a,b,\at,\bt)$ also verifies the $clos$. Finally, Lemma \ref{lem:D1eqclos} guarantees that $K(a,b,\at,\bt)$ satisfies $(D3)$.
	
	Reciprocally, assume that $K(a,b,\at,\bt)$ obeys the conditions $[D1,D2,D3]$ and is equivalent to $K$. Conditions $[D1,D2]$ ensure that $K(a,b,\at,\bt)$ verifies the $DAC$ (the details were written in the proof of Theorem \ref{teo:carDAC}). Then, thanks to Lemma \ref{lem:DACueq} we conclude that $K$ also obeys the $DAC$. Finally, the fact that $K$ verifies the $clos$ follows immediately from Lemmas \ref{lem:closueq} and \ref{lem:D1eqclos}. \qed
\end{proof}

\medskip
Since $(D1)$ and $(D3)$ provide bounds for $\at/a$ and $\bt/b$, when both of these conditions are satisfied we shortly write
\bq
\label{cond:D1p}
\mbox{$(D1')$} \qquad C_0 \le \frac{\at}{a},\ \frac{\bt}{b} \le 1.
\eq 

Clearly, $[DAC, clos] \Rightarrow [MAC, clos]$. In the following lemma we essentially prove that the converse of this estatement is also true.

\begin{lemma}
	\label{lem:red}
	Let $K$ be a convex quadrilateral obeying the $MAC$ and the $clos$, then $K$ also satisfies the $mac$. In particular, conditions $[MAC,clos]$ and $[DAC,clos]$ are equivalent.
\end{lemma}

\begin{proof}
	Let $K$ be a convex quadrilateral obeying the $MAC$ and the $clos$. Assume that $K$ has an interior angle $\beta$ tending to zero. Since the sum of all interior angles of $K$ is equal to $2\pi$, the assumption $MAC$ on $K$ implies that $\beta$ is unique. Without loss of generality, we may assume that $\beta$ is placed at $V_1$. Let $\omega_2$ and $\omega_4$ be the interior angles of $\Delta(V_1V_2V_4)$ placed at $V_2$ and $V_4$, respectively. It is clear that $\omega_2$ and $\omega_4$ are bounded away from zero and $\pi$ so, thanks to the law of sines, we conclude that $|l_{12}|$ and $|l_{14}|$ are comparable.
	Combining this fact with the assumption $clos$ it follows that $|l_{23}|$ and $|l_{34}|$ are also comparable. As consequence, all sides of $K$ have comparable lengths. From the cosine formula applied to the angle placed at $V_3$ follows that the length of the diagonal $l_{24}$ is also comparable to the length of any side of $K$. In particular, $|l_{24}|$ is comparable to $|l_{12}|$. Now, using again the law of sines on $\Delta(V_1V_2V_4)$ we derive
	$$\frac{1}{\sin(\beta)} = \frac{|l_{12}|}{|l_{24}|}\frac{1}{\sin(w_4)}$$
	and, as we shown, the right hand term is bounded above by a positive constant. This implies that $\beta$ can not tend to zero. \qed
\end{proof}

\medskip
In terms of condition $(D1')$ introduced in \eqref{cond:D1p} and thanks to Lemma \ref{lem:red}, the Theorem \ref{teo:DAC+clos} actually reads as 

\begin{theorem}
	\label{teo:MAC+clos}
	Let $K$ be a general convex quadrilateral. Then $K$ satisfies the $MAC$ and the $clos$ $($resp. $DAC$ and $clos)$ if and only if $K$ is equivalent to an element $K(a,b,\at,\bt)$ verifying conditions $[D1',D2]$ $($given by \eqref{cond:D1p} and \eqref{cond:D2}, respectively$)$.
\end{theorem}

\section{Reduction to the reference setting}
\label{sub:red}
\renewcommand{\theequation}{\arabic{section}.\arabic{equation}}
\setcounter{equation}{0}

Let $K$ be a convex quadrilateral obeying the $DAC$ and let $\oK = K(a,b,\at,\bt)$ be an equivalent element to $K$ built as was detailed in proof of Theorem \ref{teo:carDAC}. From now on, we adopt the same notation used in the proof of such theorem (in Figure \ref{fig:notonK} we partially described it).

Our next goal is to show that if the anisotropic interpolation error estimate holds on $\oK$ then a similar estimate can be obtained on $K$. In order to do this we need to guarantee that the error on $K$ is {\it comparable} to the error on $\oK$, this is the essence of the following lemma.

\begin{lemma}
\label{lem:comperror}
Let $K_1$, $K_2$ be two arbitrary convex quadrilaterals, and let $F:K_1 \to K_2$ be an affine transformation $F(x)=Bx+P$. Assume that $F(K_1)=K_2$ and $\|B\|,\|B^{-1}\| \le C$ $($i.e., $K_1$ and $K_2$ are equivalents$)$. If $Q$ and $\overline{Q}$ are the $\q_1$-interpolation on $K_1$ and $K_2$, respectively; and $\overline{u}=u \circ F^{-1}$, then for any $p \in [1,\infty)$
\bq
\label{eq:lemma22}
|\overline{u}-\overline{Q}\overline{u}|_{1,p,K_2} \lesssim 
|u-Qu|_{1,p,K_1} \lesssim 
|\overline{u}-\overline{Q}\overline{u}|_{1,p,K_2}.
\eq
\end{lemma}

\begin{proof}
We refer to \cite[Lemma 2.2]{AM} for the proof. \qed
\end{proof}

\medskip
In the sequel, we use $\overline{x}=(\ox_1,\ox_2)$ to denote the variable on $\oK$, i.e. $L(\ox_1,\ox_2)=(x_1,x_2)$ where $L : \oK=K(a,b,\at,\bt) \to K$ is the affine transformation involved in the proof of Theorem \ref{teo:carDAC} used to show that $K$ is equivalent to $\oK$. Matrix $B$ associated to $L$ is $\left( \begin{array}{cc} 1&\cot(\beta)\\ 0&1\end{array} \right)$ where $\beta$ is the interior angle of $K$ placed at $V_1$. We also use $\ov$ to denote a function defined on $\oK$ which is built starting from a function $v$ defined on $K$ by $\ov = v \circ L$. 

 On the other hand, taking into account that $l_1=V_1V_2$ and $l_2=V_1V_4$ with $V_1=(0,0)$, $V_2=(a,0)$ and $V_4=(b\cot(\beta),b)$, after a straightforward computation we obtain
\bq
\label{eq:refx1}
a \frac{\partial \ov}{\partial \ox_1} = l_1 \cdot (\nabla v \circ L) = |l_1| 
\partial_{l_1} v \circ L
\eq
and
\bq
\label{eq:refy1}
b \frac{\partial \ov}{\partial \ox_2} = l_2 \cdot (\nabla v \circ L) = |l_2| 
\partial_{l_2} v \circ L.
\eq

Changing variables and taking into account that $det(B)=1$, from 
(\ref{eq:refx1}) and (\ref{eq:refy1}) we deduce that, for any $p \in [1,\infty)$,
\bq
\label{eq:ref}
a \left\| \frac{\partial \ov}{\partial \ox_1} \right\|_{0,p,\oK} + b \left\| 
\frac{\partial \ov}{\partial \ox_2} \right\|_{0,p,\oK} = 
|l_1| \left\| \partial_{l_1} v \right\|_{0,p,K} + |l_2| \left\| \partial_{l_2} 
v \right\|_{0,p,K}.
\eq

Finally, assuming that the following anisotropic interpolation error estimate holds on $\oK$ 
\bq
\label{errorbarK}
|\ou-\oQ\ou|_{1,p,\oK} \lesssim   
a \left\| \frac{\partial}{\partial \ox_1} \nabla \ou 
\right\|_{0,p,\oK} + b \left\| \frac{\partial }{\partial \ox_2} \nabla \ou 
\right\|_{0,p,\oK},
\eq
and thanks to $L$ verifies each hypothesis in Lemma \ref{lem:comperror}, we 
get (from (\ref{eq:lemma22}) and (\ref{eq:ref}))
$$|u-Qu|_{1,p,K} \lesssim |\ou-\oQ\ou|_{1,p,\oK} \lesssim  
|l_1| \left\| \partial_{l_1} \nabla u \right\|_{0,p,K} + 
|l_2| \left\| \partial_{l_2} \nabla u \right\|_{0,p,K}.$$

This is, the following anisotropic interpolation error estimate holds on $K$
\bq
\label{eq:aee}
|u-Qu|_{1,p,K} \le C \left[ |l_1| \left\| \partial_{l_1} \nabla u 
\right\|_{0,p,K} + |l_2| \left\| \partial_{l_2} \nabla u \right\|_{0,p,K} 
\right]
\eq

where $C$ is a positive constant depending only on those constants involved in the $DAC$.

\medskip
This is the reason why the rest of this paper is devoted to prove the 
anisotropic interpolation error estimate \eqref{errorbarK} on a reference quadrilateral $\oK$ obeying conditions $[D1',D2] \Rightarrow [D1,D2]$.

\section{The error treatment}
\label{errtreat}
\setcounter{equation}{0}

In previous sections we essentially show that, in order to obtain the anisotropic interpolation error estimate on a general convex quadrilateral $K$ satisfying the $MAC$ and the $clos$, we can actually assume that $K$ belongs to the reference configuration ($K=K(a,b,\at,\bt)$) and it satisfies the conditions $[D1',D2]$. Concretely, we are interested to prove that \eqref{errorbarK} holds on $K=K(a,b,\at,\bt)$ under the assumptions $[D1',D2]$. This section is intended to present some necessary results for this purpose.

\subsection{Error decomposition}
\renewcommand{\theequation}{\arabic{section}.\arabic{subsection}.\arabic{equation}}
\setcounter{equation}{0}

Let $\Pi$ be the first-order Lagrange interpolation operator on the triangle $T_{ab}=\Delta(V_1V_2V_4)$. Then, for any $p \in [1,\infty)$, we get
\bq
\label{errordec}
|u-Qu|_{1,p,K} \le |u-\Pi u|_{1,p,K}+|\Pi u-Qu|_{1,p,K}.
\eq

Since $\Pi u - Qu$ belongs to the $\q_1$ quadrilateral finite element space and vanishes at $V_1$, $V_2$ and $V_4$, it follows that
$$\Pi u-Qu = (\Pi u - Qu)(V_3) \phi_3 = (\Pi u-u)(V_3) \phi_3$$

where $\phi_3$ is the nodal basis function associated to $V_3$. Then,
\bq
\label{destr}
|u-Qu|_{1,p,K} \le |u-\Pi u|_{1,p,K} + |(\Pi u-u)(V_3)| |\phi_3|_{1,p,K}.
\eq

\subsection{Dealing with $|(\Pi u-u)(V_3)| |\phi_3|_{1,p,K}$} 
\renewcommand{\theequation}{\arabic{section}.\arabic{subsection}.\arabic{equation}}
\setcounter{equation}{0}

In order to estimate $|(\Pi u-u)(V_3)| |\phi_3|_{1,p,K}$ we obtain bounds for $|\phi_3|_{1,p,K}$ and $|(\Pi u-u)(V_3)|$. Lemma \ref{estphi} provides an estimate to $|\phi_3|_{1,p,K}$ in terms of the diameter $h$ of $K$ and the length of $l_{34}$.

\begin{lemma}
\label{estphi}
Let $K=K(a,b,\at,\bt)$ be a convex quadrilateral satisfying $[D1,D2]$. Then, for any $p \in [1,\infty)$ with dual exponent $q$, there exists a constant $C$ depending only on $p$ and those constants involved in $[D1,D2]$ such that
$$\ds |\phi_3|_{1,p,K} \le C \frac{h^{1/p}}{|l_{34}|^{1/q}}.$$
\end{lemma}

\begin{proof}
See \cite[Lemma 5.3 (d)]{AM}. \qed
\end{proof}

\medskip
The treatment of $|(\Pi u-u)(V_3)|$ is slightly more difficult and requires the use of a trace theorem that we give in Theorem \ref{teo:trace}. Such theorem can be regarded as a generalized sharp version of \cite[Lemma 3.2]{Verf} and its proof relies on the following

\begin{lemma}
\label{lem:verfaux}
Let $\hT$ be the $2$-simplex, i.e. $\hT=\Delta(\widehat{V}^1\widehat{V}^2\widehat{V}^3)$ where $\widehat{V}^1=(0,0)$, $\widehat{V}^2=(1,0)$ and $\widehat{V}^3=(0,1)$. For any $p \in [1,\infty)$ and $v \in W^{1,p}(\hT)$ which vanishes on $\widehat{V}^2\widehat{V}^3$ we have
$$\left\| v \right\|_{0,p,l_1} \le \left\| \frac{\partial v}{\partial x_2} \right\|_{0,p,\hT} \qquad \text{and} \qquad 
\left\| v \right\|_{0,p,l_2} \le \left\| \frac{\partial v}{\partial x_1} \right\|_{0,p,\hT}$$
where $l_1=\widehat{V}^1\widehat{V}^2$ and $l_2=\widehat{V}^1\widehat{V}^3$.
\end{lemma}

\begin{proof}
This result is a generalization of Lemma 3.1 in \cite{Verf} to any $p \in [1,\infty)$ when $n=2$ and its proof is a straightforward adaptation of the one given in \cite{Verf} so we omit it. \qed
\end{proof}

\begin{theorem}
	\label{teo:trace}	
	Let $T$ be the triangle of vertices $V^1, V^2$ and $V^3$ where $V^i = (V^i_1, V^i_2)$, $i=1, 2, 3$. Assume the existence of positive constants $a_1, a_2$ such that 
	\bq
	\label{aj}
	|V^i_j-V^k_j| \le a_j, \quad 1 \le k < i \le 3, \quad j=1,2.
	\eq	
	
	If $\lambda_i$ is the barycentric coordinate associated with $V^i$, then for any $p \in [1,\infty)$ and $v \in W^{1,p}(T)$ 
	\bq
	\label{eq:vi}
	\left\| \lambda_i v \right\|_{0,p,l} \le
	2^{\frac{p-1}{p}}\left( \frac{|l|}{|T|} \right)^{1/p} \left[ \left\| v \right\|_{0,p,T}+
	a_1 \left\| \frac{\partial v}{\partial x_1} \right\|_{0,p,T} + 
	a_2 \left\| \frac{\partial v}{\partial x_2} \right\|_{0,p,T} \right]
	\eq
	where $l$ is any edge of $T$ having to $V^i$ as an extreme.
	
	In addition, if $e$ is any edge of $T$, then for any $p \in [1,\infty)$ and $v \in W^{1,p}(T)$ 
	\bq
	\label{eq:e}
	\left\| v \right\|_{0,p,e} \le
	2^{\frac{2p-1}{p}}\left( \frac{|l|}{|T|} \right)^{1/p} \left[ \left\| v \right\|_{0,p,T}+
	a_1 \left\| \frac{\partial v}{\partial x_1} \right\|_{0,p,T} + 
	a_2 \left\| \frac{\partial v}{\partial x_2} \right\|_{0,p,T} \right].
	\eq
\end{theorem}

\begin{proof}
	We only prove \eqref{eq:vi} for $i=1$ since the arguments can be easily adapted to obtain the result for $i=2,3$ changing appropriately the mapping $F_T$ given in \eqref{def:FT} below.
	
	Let $\hT$ the triangle of vertices $\widehat{V}^1=(0,0)$, $\widehat{V}^2=(1,0)$ and $\widehat{V}^3=(0,1)$. Consider the affine transformation which maps $\hT$ onto $T$ and $\widehat{V}^1$ to $V^1$, $F_T:\hT \to T$, given by
	\bq
	\label{def:FT}
	F_T(\hx) = \begin{pmatrix} V^2-V^1\ |\ V^3-V^1 \end{pmatrix} \hx + V^1.
	\eq
	
	Let $l = V^1V^{j+1}$ with $j=1$ or $j=2$. Then $l$ is the image w.r.t. $F_T$ of segment $\hat{l} = \widehat{V}^1 \widehat{V}^{j+1}$; on the other hand, $\lambda_1 \circ F_T = \widehat{\lambda}_1$ where $\widehat{\lambda}_1$ is the barycentric coordinate of $\hT$ associated with $\widehat{V}^1$. Set $\widehat{v}=v \circ F_T$. Let $i=3-j$, since $\widehat{\lambda}_1 \widehat{v}$ vanishes on the edge $\widehat{V}^2\widehat{V}^3$ of $\hT$ we may apply Lemma \ref{lem:verfaux} and obtain
	$$\begin{array}{lcl}
	\left\| \lambda_1 v \right\|_{0,p,l} &=& \ds \left( \frac{|l|}{|\hat{l}|} \right)^{1/p} \left\| \widehat{\lambda}_1 \widehat{v} \right\|_{0,p,\hat{l}} \medskip\\
	&\le& \ds \left( \frac{|l|}{|\hat{l}|} \right)^{1/p} \left\| \frac{\partial}{\partial \hx_i} (\widehat{\lambda}_1 \widehat{v}) \right\|_{0,p,\hT} \medskip\\  
	&\le& \ds \left( \frac{|l|}{|\hat{l}|} \right)^{1/p} \left[ \left\| \widehat{v} \frac{\partial \widehat{\lambda}_1}{\partial \hx_i} \right\|_{0,p,\hT} + 
	\left\| \widehat{\lambda}_1 \frac{\partial \widehat{v}}{\partial \hx_i} \right\|_{0,p,\hT} \right].
	\end{array}$$
	
	Since $\ds \frac{\partial \widehat{\lambda}_1}{\partial \hx_i}=-1$ and $\left\| \widehat{\lambda}_1 \right\|_{\infty,\hT}=1$ this yields
	\bq
	\label{est:lp1}
	\left\| \lambda_1 v \right\|_{0,p,l} \le \left( \frac{|l|}{|\hat{l}|} \right)^{1/p} \left[ \left\| \widehat{v} \right\|_{0,p,\hT} + 
	\left\| \frac{\partial \widehat{v}}{\partial \hx_i} \right\|_{0,p,\hT} \right].
	\eq
	
	Transforming back to $\hT$ we get
	\bq
	\label{eq:vtl2}
	\left\| \widehat{v} \right\|_{0,p,\hT} = 
	\left( \frac{|\hT|}{|T|} \right)^{1/p} \left\| v \right\|_{0,p,T}
	\eq
	and 
	$$\left\| \frac{\partial \widehat{v}}{\partial \hx_i} \right\|_{0,p,\hT} =
	\left( \frac{|\hT|}{|T|} \right)^{1/p} \left\| \nabla v \cdot DF_T \cdot E_i^T \right\|_{0,p,T}$$
	where $E_i$ denotes the $i$th unit vector. 
	
	Now, taking into account that $DF_T \cdot E_i^T = \begin{pmatrix} V^2-V^1\ |\ V^3-V^1 \end{pmatrix} \cdot E_i^T = V^{i+1}-V^1$ and using the $p$-inequality $(s+t)^p \le 2^{p-1}(s^p+t^p)$ ($s,t \ge 0$ and $p \in [1,\infty)$) we get
	$$\begin{array}{lcl}
	\left\| \nabla v \cdot DF_T \cdot E_i^T \right\|_{0,p,T}^p &=& 
	\left\| \nabla v \cdot \begin{pmatrix} V^{i+1}-V^1 \end{pmatrix} \right\|_{0,p,T}^p \medskip\\ &\le& \ds 
	2^{p-1} \left[ \left| V^{i+1}_1-V^1_1 \right|^p \left\| \frac{\partial v}{\partial x_1} \right\|_{0,p,T}^p + 
	\left| V^{i+1}_2-V^1_2 \right|^p \left\| \frac{\partial v}{\partial x_2} \right\|_{0,p,T}^p \right].
	\end{array}$$
	
	Then, from the assumption \eqref{aj} we obtain
	$$\ds \left\| \nabla v \cdot DF_T \cdot E_i^T \right\|_{0,p,T}^p \le  
	2^{p-1} \left[ a_1^p \left\| \frac{\partial v}{\partial x_1} \right\|_{0,p,T}^p + 
	a_2^p \left\| \frac{\partial v}{\partial x_2} \right\|_{0,p,T}^p \right],$$
	and hence 
	\bq
	\label{eq:derl2}
	\left\| \frac{\partial \widehat{v}}{\partial \hx_i} \right\|_{0,p,\hT} \le
	2 \left( \frac{|\hT|}{|T|} \right)^{1/p} \left[ 
	a_1 \left\| \frac{\partial v}{\partial x_1} \right\|_{0,p,T} + 
	a_2 \left\| \frac{\partial v}{\partial x_2} \right\|_{0,p,T} \right].
	\eq
	
	Taking into account that $\ds \left( \frac{|\hT|}{|\hat{l}|} \right)^{1/p}=2^{-1/p}$, and combining \eqref{est:lp1} with \eqref{eq:vtl2} and \eqref{eq:derl2} we have
	$$\left\| \lambda_1 v \right\|_{0,p,l} \le
	2^{\frac{p-1}{p}} \left( \frac{|l|}{|T|} \right)^{1/p} 
	\left[ \left\| v \right\|_{0,p,T} + 
	a_1 \left\| \frac{\partial v}{\partial x_1} \right\|_{0,p,T} + 
	a_2 \left\| \frac{\partial v}{\partial x_2} \right\|_{0,p,T} \right]$$
	which concludes the proof of \eqref{eq:vi} for $i=1$.
	
	Finally, \eqref{eq:e} follows immediately from \eqref{eq:vi}. Indeed, let $e$ be any edge of $T$; renaming the vertices if necessary we can assume $e=V^1V^2$. Taking into account that any function $v$ can be written as $v=(\lambda_1+\lambda_2+\lambda_3)v$ and $\lambda_3$ vanishes on $e$ we have, on $e$, $v=(\lambda_1+\lambda_2)v$. Then, using triangular inequality it follows that $\left\| v \right\|_{0,p,e} \le \left\| \lambda_1 v \right\|_{0,p,e}+\left\| \lambda_2 v \right\|_{0,p,e}$; therefore, \eqref{eq:e} follows by using the estimates given in \eqref{eq:vi} for $i=1, 2$. \qed
	\end{proof} 

\begin{figure}[h]
	\centering
	\resizebox{16cm}{7cm}{\includegraphics{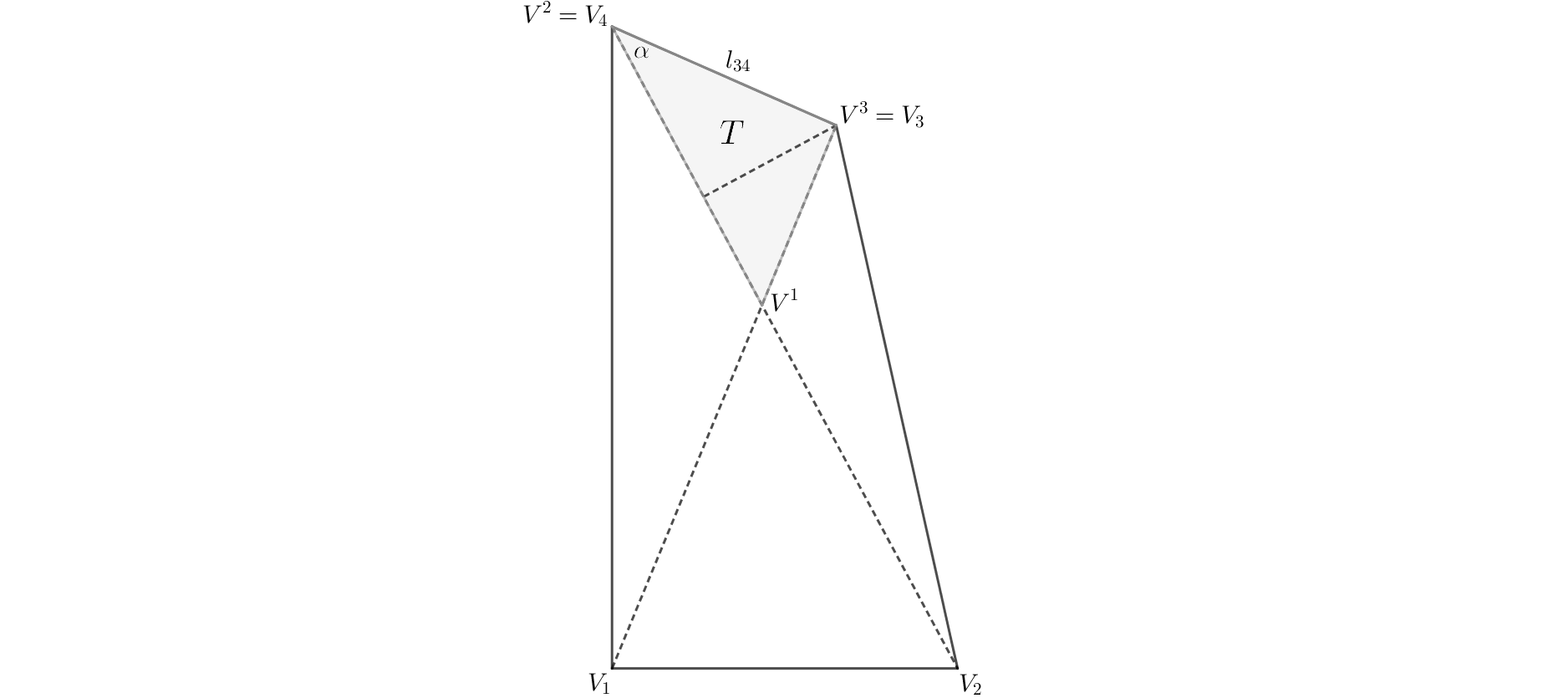}}	
	\caption{Gray area denotes the triangle $T=\Delta(V^1V^2V^3)$ where $V^1=l_{13} \cap l_{24}$, $V^2=V_4$ and $V^3=V_3$.}
	\label{fig:traceT}
\end{figure}
\medskip
We use Theorem \ref{teo:trace} to bound $|(u-\Pi u)(V_3)|$. The corresponding estimate is estated in the following lemma.

\begin{lemma}
\label{upikumij}
Let $K=K(a,b,\at,\bt)$ be a convex quadrilateral satisfying conditions $[D1',D2]$. For any $p \in [1,\infty)$ with dual exponent $q$, we have 
$$|(u-\Pi u)(V_3)| \le \left( \frac{2^{3p+1}C}{C_0} \right)^{1/p} \frac{|l_{34}|^{1/q}}{h^{1/p}} \left[ \|\partial_{l_{34}}(u-\Pi u)\|_{0,p,T} + a \left\| \frac{\partial}{\partial x_1} \nabla u \right\|_{0,p,T} + b \left\| \frac{\partial}{\partial x_2} \nabla u \right\|_{0,p,T} \right]$$
where $T=\Delta(V^1V^2V^3)$ with $V^1=l_{13} \cap l_{24}$, $V^2=V_4$, $V^3=V_3$ $($see {\rm Figure \ref{fig:traceT})} and $C_0,C$ are the constants involved in $(D1')$ and $(D2)$, respectively.
\end{lemma}

\begin{proof}
In order to apply Theorem \ref{teo:trace} on $T=\Delta(V^1V^2V^3)$ we need to guarantee the imposed requirement \eqref{aj}. A simple computation shows that $V^1=\left( a \frac{\at b}{\at b+a\bt},  b \frac{a\bt}{\at b+a\bt} \right)$; then, taking into account that $0<\frac{\at b}{\at b+a\bt}, \frac{a \bt}{\at b+a\bt} \le 1$ we conclude $|V_1^1| \le a$ and $|V_2^1| \le b$. On the other hand, since $V^2 = (0,b)$ and $V^3=(\at,\bt)$, thanks to $(D1')$, we also have $|V^j_1| \le a$ and $|V^j_2| \le b$, $j=2,3$. Finally, for any $1 \le k \le i \le 3$, we obtain (by using triangular inequality)
$$|V^i_1-V^k_1| \le 2a=:a_1 \qquad \mbox{and} \qquad 
|V^i_2-V^k_2| \le 2b=:a_2,$$
and hence the condition \eqref{aj} is fulfilled.

\medskip	
Now, since $(u-\Pi u)(V_4)=0$, a combination of the H\"older's inequality and the Theorem \ref{teo:trace} yield
$$\begin{array}{lcl}
|(u - \Pi u)(V_3)| &=& \ds \Big| \int_{l_{34}} \partial_{l_{34}}(u-\Pi u) \Big| \medskip\\
&\le& |l_{34}|^{1/q} \|\partial_{l_{34}}(u-\Pi u)\|_{0,p,l_{34}} \medskip\\
&\le& \ds 2^{\frac{3p-1}{p}} \frac{|l_{34}|}{|T|^{1/p}} 
\left[ \|\partial_{l_{34}}(u-\Pi u)\|_{0,p,T} + a \left\| \frac{\partial}{\partial x_1} \nabla u \right\|_{0,p,T} + b \left\| \frac{\partial}{\partial x_2} \nabla u \right\|_{0,p,T} \right].
\end{array}$$

On the other hand, from \eqref{heqv2v4} it follows that
\bq
\label{v1v2}
\ds |V^1V^2| = \frac{\at b}{\at b+a\bt} (a^2+b^2)^{1/2} = \frac{1}{1+\frac{a}{\at} \frac{\bt}{b}} h.
\eq

Condition $(D1')$ implies that $1+\frac{a}{\at} \frac{\bt}{b} \le 1+\frac{a}{\at} \le \frac{2}{C_0}$, which combined with \eqref{v1v2} give us
\bq
\label{v1v2final}
\ds |V^1V^2| \ge \frac{C_0}{2} h.
\eq

Writing $|T| = \frac{1}{2}|V^1V^2| |l_{34}|\sin(\alpha)$ (see Figure \ref{fig:traceT}), we obtain $|T| \ge \frac{C_0}{4} h |l_{34}| \sin(\alpha)$. Then, thanks to $(D2)$, it follows that 
$$\frac{|l_{34}|}{|T|^{1/p}} \le 
\left( \frac{2^2}{C_0} \right)^{1/p} \frac{|l_{34}|^{1/q}}{(\sin(\alpha)h)^{1/p}} \le
\left( \frac{2^2C}{C_0} \right)^{1/p} \frac{|l_{34}|^{1/q}}{h^{1/p}}$$
which completes the proof. \qed
\end{proof}

\medskip
Finally, Lemmas \ref{estphi} and \ref{upikumij} allow to estimate $|(u-\Pi u)(V_3)| |\phi_3|_{1,p,K}$.

\begin{lemma}
	\label{lem:upiuphi3}
	Let $K=K(a,b,\at,\bt)$ be a convex quadrilateral satisfying conditions $[D1',D2]$. For any $p \in [1,\infty)$ with dual exponent $q$, there exists a constant $C$ such that 
	$$|(u-\Pi u)(V_3)| |\phi_3|_{1,p,K} \le C \left[ \|\partial_{l_{34}}(u-\Pi u)\|_{0,p,T} + a \left\| \frac{\partial}{\partial x_1} \nabla u \right\|_{0,p,T} + b \left\| \frac{\partial}{\partial x_2} \nabla u \right\|_{0,p,T} \right]$$
	where $T=\Delta(V^1V^2V^3)$ with $V^1=l_{13} \cap l_{24}$, $V^2=V_4$ and $V^3=V_3$.
\end{lemma}

\subsection{Dealing with $|u-\Pi u|_{1,p,K}$}
\setcounter{equation}{0}

Since $K$ is the union of $T_{234}=\Delta(V_2V_3V_4)$ and $T_{ab}=\Delta(V_1V_2V_4)$ (with $|T_{234}\cap T_{ab}|=0$) it follows that
\bq
\label{kasunion}
|u-\Pi u|_{1,p,K} = |u-\Pi u|_{1,p,T_{ab}}+|u-\Pi u|_{1,p,T_{234}}.
\eq

Let $\Pi^{T_{234}}$ be the first-order Lagrange interpolation operator on $T_{234}$, then (by using triangular inequality at the second term on the right hand side) we obtain
\bq
\label{decoerroronT}
|u-\Pi u|_{1,p,K} \le |u-\Pi u|_{1,p,T_{ab}}+|u-\Pi^{T_{234}} u|_{1,p,T_{234}}+|\Pi^{T_{234}} u-\Pi u|_{1,p,T_{234}}.
\eq

As we shall see soon, terms $|u-\Pi u|_{1,p,T_{ab}}$ and $|u-\Pi^{T_{234}} u|_{1,p,T_{234}}$ can be easily bounded thanks to known results on triangles; this is why our focus is on $|\Pi^{T_{234}} u-\Pi u|_{1,p,T_{234}}$.

Notice that function $\Pi^{T_{234}} u-\Pi u$ belongs to $\mathbb{P}_1$ so it can be written as a linear combination of $\{ \lambda_2, \lambda_3, \lambda_4 \}$ where $\lambda_i$ is the barycentric coordinate on $T_{234}=\Delta(V_2V_3V_4)$ associated to $V_i$. Taking into account that $\Pi^{T_{234}} u$ agrees with $\Pi u$ on $V_2$ and $V_4$ we deduce
$$\Pi^{T_{234}} u-\Pi u = (\Pi^{T_{234}} u-\Pi u)(V_3) \lambda_3 = (u-\Pi u)(V_3) \lambda_3.$$

Then
\bq
\label{lambda3sem1}
|\Pi^{T_{234}} u-\Pi u|_{1,p,T_{234}} \le |(u-\Pi u)(V_3)| |\lambda_3|_{1,p,T_{234}}.
\eq

\medskip
Lemma \ref{upikumij} provides an estimate for $|(u-\Pi u)(V_3)|$; on the other hand, w.r.t $|\lambda_3|_{1,p,T_{234}}$ we can bound this term in a similar fashion than $|\phi_3|_{1,p,K}$ was bounded in Lemma \ref{estphi}. 

\begin{lemma}
\label{estlambda}
Let $K=K(a,b,\at,\bt)$ be a convex quadrilateral satisfying $[D1,D2]$. Then, for any $p \in [1,\infty)$ with dual exponent $q$, we have
$$\ds |\lambda_3|_{1,p,T_{234}} \le C \frac{h^{1/p}}{|l_{34}|^{1/q}}$$
where $C$ only depends on $p$ and those constants involved in $[D1,D2]$.	
\end{lemma}

\begin{proof}
A simple computation shows that $\lambda_3(x)=1+\frac{b(x_1-\at)+a(x_2-\bt)}{JF_{T_{234}}}$ where $JF_{T_{234}}$ is the Jacobian of the affine mapping $F_{T_{234}}:\hT \to T_{234}$ defined by $F_{T_{234}}(\hx) = \begin{pmatrix} V_3-V_2\ |\ V_4-V_2 \end{pmatrix} \hx + V_2$.

Since $|JF_{T_{234}}|=\dfrac{|T_{234}|}{|\hT|}=2|T_{234}|$ we have
\bq
\label{lambda3}
|\lambda_3|_{1,p,T_{234}} = \frac{|T_{234}|^{1/p}}{|JF_{T_{234}}|} (a+b) = 2^{-1} \frac{a+b}{|T_{234}|^{1/q}} \le \frac{h}{|T_{234}|^{1/q}}.
\eq

From \eqref{heqv2v4} we know that $h=|l_{24}|$ and then $|T_{234}| = \frac{1}{2} |l_{34}| h \sin(\alpha)$. Combining this fact with \eqref{lambda3}, and thanks to $(D2)$, we conclude
\bq
\label{boundlambda3}
|\lambda_3|_{1,p,T_{234}} \le 2^{1/q}\frac{h}{(|l_{34}|h\sin(\alpha))^{1/q}} \le 
(2C)^{1/q}\frac{h^{1/p}}{|l_{34}|^{1/q}}
\eq
where $C$ is the constant involved in $(D2)$. \qed
\end{proof}

\medskip
From \eqref{lambda3sem1} and Lemmas \ref{upikumij} and \ref{estlambda}, we obtain the following estimate for $|\Pi^{T_{234}} u-\Pi u|_{1,p,T_{234}}$.

\begin{lemma}
\label{lem:PiTPi}
Let $K=K(a,b,\at,\bt)$ be a convex quadrilateral satisfying $[D1',D2]$. For any $p \in [1,\infty)$, there exists a positive constant $C$ such that  
$$|\Pi^{T_{234}} u-\Pi u|_{1,p,T_{234}} \le C \left[ \|\partial_{l_{34}}(u-\Pi u)\|_{0,p,T} + a \left\| \frac{\partial}{\partial x_1} \nabla u \right\|_{0,p,T} + b \left\| \frac{\partial}{\partial x_2} \nabla u \right\|_{0,p,T} \right],$$
where $T=\Delta(V^1V^2V^3)$ with $V^1=l_{13} \cap l_{24}$, $V^2=V_4$ and $V^3=V_3$.
\end{lemma}

\section{Anisotropic error estimates on triangles}
\label{ontriangles}
\renewcommand{\theequation}{\arabic{section}.\arabic{equation}}
\setcounter{equation}{0}

From \eqref{destr}, \eqref{decoerroronT}, Lemma \ref{lem:upiuphi3} and Lemma \ref{lem:PiTPi} we deduce that the interpolation error estimate on $K=K(a,b,\at,\bt)$ reduces to obtain suitable bounds for $|u-\Pi u|_{1,p,T_{ab}}$, $|u-\Pi^{T_{234}} u|_{1,p,T_{234}}$ and $\| \partial_{l_{34}} (u-\Pi u) \|_{0,p,T}$. The treatment of these terms relies on Lemma \ref{lem:ric} which is essentially a Poincar\'e type inequality on triangles.

\begin{lemma}
\label{lem:ric}
Let $T$ be a triangle with edges $e$, $l_1$ and $l_2$. Let $v \in W^{1,p}(T)$, $p \in [1,\infty)$, be a function with vanishing average on $e$. Then, there exists a constant $C$ independent of $T$ such that
$$\left\| v \right\|_{0,p,T} \le C \left[ 
|l_1| \left\| \partial_{l_1} v \right\|_{0,p,T} +  
|l_2| \left\| \partial_{l_2} v \right\|_{0,p,T} \right].$$
\end{lemma}

\begin{proof}
For $p=2$, this result reduces to \cite[Lemma 2.2]{D} and the proof given there can be replied step by step to any $p$. \qed
\end{proof}

\begin{lemma}
\label{lem:uPiuonT}
Let $K=K(a,b,\at,\bt)$ be a convex quadrilateral satisfying $[D1,D2]$ and let $T$ be the triangle $\Delta(V^1V^2V^3)$ where $V^1=l_{13} \cap l_{24}$, $V^2=V_4$, $V^3=V_3$. Then, for any $p \in [1,\infty)$, there exists a positive constant $C$ $($depending only on $p$ and those constants involved in $(D1)$ and $(D2))$ such that  
$$\left\| \partial_{l_{34}} (u-\Pi u) \right\|_{0,p,T} \le C \left[ 
a \left\| \frac{\partial}{\partial x_1} \nabla u \right\|_{0,p,K} +
b \left\| \frac{\partial}{\partial x_2} \nabla u \right\|_{0,p,K} \right].$$
\end{lemma}

\begin{proof}
Notice that $\partial_{l_{14}} (u-\Pi u)$ has vanishing average on $l_{14}$ due to $(u-\Pi u)(V_1)=0=(u-\Pi u)(V_4)$; then we may apply Lemma \ref{lem:ric} on $T_{134}=\Delta(V_1V_3V_4)$ and, using that the second derivatives of $\Pi u$ vanish, we obtain
$$\left\| \partial_{l_{14}} (u-\Pi u) \right\|_{0,p,T_{134}} \lesssim 
|l_{13}| \left\|\partial_{l_{13}} \partial_{l_{14}} u \right\|_{0,p,T_{134}} +
|l_{34}| \left\|\partial_{l_{34}} \partial_{l_{14}} u \right\|_{0,p,T_{134}}.$$

Since $V_1=(0,0)$ and $V_3=(\at,\bt)$, by using the triangular inequality and the fact $(s+t)^p \le 2^{p-1}(s^p+t^p)$ ($s,t \ge 0$ and $p \in [1,\infty)$), we deduce

$$\begin{array}{lcl}
|l_{13}| \left\|\partial_{l_{13}} \partial_{l_{14}} u \right\|_{0,p,T_{134}} & \le &
\ds \at \left\| \frac{\partial}{ \partial x_1}(\partial_{l_{14}} u) \right\|_{0,p,T_{134}} + 
\bt \left\| \frac{\partial}{ \partial x_2}(\partial_{l_{14}} u) \right\|_{0,p,T_{134}} \medskip\\
& = &
\ds \at \left\| \frac{\partial}{ \partial x_1}(\nabla u) \cdot \frac{l_{14}}{|l_{14}|} \right\|_{0,p,T_{134}} + 
\bt \left\| \frac{\partial}{ \partial x_2}(\nabla u) \cdot \frac{l_{14}}{|l_{14}|} \right\|_{0,p,T_{134}} \medskip\\
& \lesssim &
\ds \at \left\| \frac{\partial}{ \partial x_1}(\nabla u) \right\|_{0,p,T_{134}} + 
\bt \left\| \frac{\partial}{ \partial x_2}(\nabla u) \right\|_{0,p,T_{134}}.
\end{array}$$

Finally, from $(D1)$ it follows that 
$$|l_{13}| \left\|\partial_{l_{13}} \partial_{l_{14}} u \right\|_{0,p,T_{134}} \lesssim a \left\| \frac{\partial}{ \partial x_1}(\nabla u) \right\|_{0,p,T_{134}} + b \left\| \frac{\partial}{ \partial x_2}(\nabla u) \right\|_{0,p,T_{134}}.$$

In a similar fashion we have $|l_{34}| \left\|\partial_{l_{34}} \partial_{l_{14}} u \right\|_{0,p,T_{134}} \lesssim a \left\| \frac{\partial}{ \partial x_1}(\nabla u) \right\|_{0,p,T_{134}} + b \left\| \frac{\partial}{ \partial x_2}(\nabla u) \right\|_{0,p,T_{134}}$. Hence,
\bq
\label{dx2upiu2}
\left\| \partial_{l_{14}} (u-\Pi u) \right\|_{0,p,T_{134}} \lesssim 
a \left\| \frac{\partial}{\partial x_1} \nabla u \right\|_{0,p,T_{134}} +
b \left\| \frac{\partial}{\partial x_2} \nabla u \right\|_{0,p,T_{134}}.
\eq

On the other hand, the function $\partial_{l_{24}} (u-\Pi u)$ has vanishing average on $l_{24}$ since $(u-\Pi u)(V_2)=0=(u-\Pi u)(V_4)$; then we may apply Lemma \ref{lem:ric} on triangle  $T_{234}=\Delta(V_2V_3V_4)$ to obtain
$$\left\| \partial_{l_{24}} (u-\Pi u) \right\|_{0,p,T_{234}} \lesssim 
|l_{23}| \left\|\partial_{l_{23}} \partial_{l_{24}} u \right\|_{0,p,T_{234}} +
|l_{34}| \left\|\partial_{l_{34}} \partial_{l_{24}} u \right\|_{0,p,T_{234}}.$$

Proceeding as before we obtain
$$|l_{23}| \left\|\partial_{l_{23}} \partial_{l_{24}} u \right\|_{0,p,T_{234}}, 
|l_{34}| \left\|\partial_{l_{34}} \partial_{l_{24}} u \right\|_{0,p,T_{234}}
\lesssim a \left\| \frac{\partial}{ \partial x_1}(\nabla u) \right\|_{0,p,T_{234}} + b \left\| \frac{\partial}{ \partial x_2}(\nabla u) \right\|_{0,p,T_{234}},$$
and hence,
\bq
\label{dl24upiu2}
\left\| \partial_{l_{24}} (u-\Pi u) \right\|_{0,p,T_{234}} \lesssim 
a \left\| \frac{\partial}{\partial x_1} \nabla u \right\|_{0,p,T_{234}} +
b \left\| \frac{\partial}{\partial x_2} \nabla u \right\|_{0,p,T_{234}}.
\eq

From inclusions $T \subset T_{i34} \subset K$ for $i=1, 2$, and \eqref{dx2upiu2}, \eqref{dl24upiu2} we derive
\bq
\label{dli4upiuT}
\left\| \partial_{l_{i4}} (u-\Pi u) \right\|_{0,p,T} \lesssim 
a \left\| \frac{\partial}{\partial x_1} \nabla u \right\|_{0,p,K} +
b \left\| \frac{\partial}{\partial x_2} \nabla u \right\|_{0,p,K} \qquad i=1, 2.
\eq

Notice that $l_{34} = (1-\at/a-\bt/b) l_{14} + \at/a\ l_{24}$. Then, using triangular inequality and $(D1)$, it follows that $\left| \partial_{l_{34}} (u-\Pi u) \right| \lesssim 
\left| \partial_{l_{14}} (u-\Pi u) \right| + 
\left| \partial_{l_{24}} (u-\Pi u) \right|$.

Finally, using one more time the fact $(s+t)^p \le 2^{p-1}(s^p+t^p)$ ($s,t \ge 0$ and $p \in [1,\infty)$), from \eqref{dli4upiuT} we deduce
$$\left\| \partial_{l_{34}} (u-\Pi u) \right\|_{0,p,T} \lesssim  
a \left\| \frac{\partial}{\partial x_1} \nabla u \right\|_{0,p,K} +
b \left\| \frac{\partial}{\partial x_2} \nabla u \right\|_{0,p,K}.$$ \qed 
\end{proof}

\begin{lemma}
\label{lem:errorontriangles}
Let $K=K(a,b,\at,\bt)$ be a convex quadrilateral satisfying conditions $[D1,D2]$. Let $T_{ab}$ and $T_{234}$ be the triangles $\Delta(V_1V_2V_4)$ and $\Delta(V_2V_3V_4)$, respectively. Then, for any $p \in [1,\infty)$, there exists a positive constant $C$ $($depending only on $p$ and those constants involved in $(D1)$ and $(D2))$ such that  
\bq
\label{errorsobreTab}
| u-\Pi u |_{1,p,T_{ab}} \le C \left[ 
a \left\| \frac{\partial}{\partial x_1} \nabla u \right\|_{0,p,K} +
b \left\| \frac{\partial}{\partial x_2} \nabla u \right\|_{0,p,K} \right]
\eq
and
\bq
\label{errorsobreT234}
| u-\Pi^{T_{234}} u |_{1,p,T_{234}} \le C \left[ 
a \left\| \frac{\partial}{\partial x_1} \nabla u \right\|_{0,p,K} +
b \left\| \frac{\partial}{\partial x_2} \nabla u \right\|_{0,p,K} \right].
\eq
\end{lemma}

\begin{proof}
Since $\partial_{l_{24}}(u-\Pi^{T_{234}} u)$ and $\partial_{l_{34}}(u-\Pi^{T_{234}} u)$ have vanishing average on $l_{24}$ and $l_{34}$, respectively (due to $(u-\Pi^{T_{234}} u)(V_i)=0$ for $i=2, 3, 4$), we may apply Lemma \ref{lem:ric} and, using that the second derivatives of $\Pi^{T_{234}} u$ vanish, we obtain
$$\left\| \partial_{l_{24}} (u-\Pi^{T_{234}} u) \right\|_{0,p,T_{234}} \lesssim 
|l_{23}| \left\|\partial_{l_{23}} \partial_{l_{24}} u \right\|_{0,p,T_{234}} +
|l_{34}| \left\|\partial_{l_{34}} \partial_{l_{24}} u \right\|_{0,p,T_{234}}$$
and
$$\left\| \partial_{l_{34}} (u-\Pi^{T_{234}} u) \right\|_{0,p,T_{234}} \lesssim 
|l_{23}| \left\|\partial_{l_{23}} \partial_{l_{34}} u \right\|_{0,p,T_{234}} +
|l_{24}| \left\|\partial_{l_{24}} \partial_{l_{34}} u \right\|_{0,p,T_{234}}.$$

With similar arguments than those used in the proof of Lemma \ref{lem:uPiuonT}, for any term on the right in previous inequalities, we conclude that
$$|l_{ij}| \left\|\partial_{l_{ij}} \partial_{l_{kr}} u \right\|_{0,p,T_{234}} \lesssim a \left\| \frac{\partial}{ \partial x_1}(\nabla u) \right\|_{0,p,T_{234}} + b \left\| \frac{\partial}{ \partial x_2}(\nabla u) \right\|_{0,p,T_{234}}.$$

Then, for $i=2, 3$, we have
\bq
\label{comb}
\left\| \partial_{l_{i4}} (u-\Pi^{T_{234}} u) \right\|_{0,p,T_{234}} \lesssim 
a \left\| \frac{\partial}{\partial x_1} \nabla u \right\|_{0,p,T_{234}} +
b \left\| \frac{\partial}{\partial x_2} \nabla u \right\|_{0,p,T_{234}}.
\eq

Finally, the estimate \eqref{errorsobreT234} for $|u-\Pi^{T_{234}}|_{1,p,T_{234}}$ follows from \eqref{comb} by using that the angle $\alpha$ between $l_{24}$ and $l_{34}$ is bounded away from zero and $\pi$ thanks to $(D2)$; and the fact that $T_{234}$ is contained in $K$. 

Estimate \eqref{errorsobreTab} can be obtained in a similar way. Indeed, notice that $\partial_{l_{24}}(u-\Pi u)$ and $\partial_{l_{12}}(u-\Pi u)$ have vanishing average on $l_{24}$ and $l_{12}$, respectively. Then we may apply Lemma \ref{lem:ric} and, bounding the corresponding terms as we did before, we obtain
\bq
\label{comb2}
\left\| \partial_{l_{24}} (u-\Pi u) \right\|_{0,p,T}, \left\| \partial_{l_{12}} (u-\Pi u) \right\|_{0,p,T} \lesssim 
a \left\| \frac{\partial}{\partial x_1} \nabla u \right\|_{0,p,T} +
b \left\| \frac{\partial}{\partial x_2} \nabla u \right\|_{0,p,T}.
\eq

The proof concludes by noticing that the angle between $l_{12}$ and $l_{24}$ is greater or equal than $\alpha$ so is bounded away from zero thanks to $(D2)$, and is lower than $\pi/2$ thanks to $(D1)$; together with the elementary fact that $T_{ab}$ is contained in $K$. \qed
\end{proof}

\section{Main results}
\label{main}
\setcounter{equation}{0}

In this section we present our main result (Theorem \ref{teo:main}) which deals with arbitrary convex quadrilaterals satisfying the $MAC$ and the $clos$; however, in order to keep our exposure as clear as possible, we begin by proving an auxiliary result on a subclass of elements belonging to the reference setting which is interesting by itself and is written in such a way to attend our purposes.

\begin{theorem}
\label{teo:ref}
Let $K=K(a,b,\at,\bt)$ be a convex quadrilateral satisfying conditions $[D1',D2]$. For any $p \in [1,\infty)$, there exists a constant $C$ $($depending only on $p$ and those constants involved in $(D1')$ and $(D2))$ such that
\bq
\label{eq:teoref}
|u-Qu|_{1,p,K} \le C \left[ 
a \left\| \frac{\partial}{\partial x_1} \nabla u \right\|_{0,p,K} + 
b \left\| \frac{\partial}{\partial x_2} \nabla u \right\|_{0,p,K} \right].
\eq
\end{theorem}

\begin{proof}
As was detailed in \eqref{destr} and \eqref{decoerroronT}, the term $|u-Qu|_{1,p,K}$ can be descomposed as follows
$$\begin{array}{lcl}
|u-Qu|_{1,p,K} &\le& |u-\Pi u|_{1,p,K} + |(\Pi u-u)(V_3)| |\phi_3|_{1,p,K} \medskip\\
&\le& |u-\Pi u|_{1,p,T_{ab}}+|u-\Pi^{T_{234}} u|_{1,p,T_{234}}+|\Pi^{T_{234}} u-\Pi u|_{1,p,T_{234}}+ \medskip\\
&  & \hspace{8.5cm} +|(\Pi u-u)(V_3)| |\phi_3|_{1,p,K}
\end{array}$$
where $\Pi$ is the first-order Lagrange interpolation operator on $T_{ab}=\Delta(V_1V_2V_4)$ and $\Pi^{T_{234}}$ is the first-order Lagrange interpolation operator on $T_{234}=\Delta(V_2V_3V_4)$.

Taking into account that conditions $[D1',D2]$ imply $[D1,D2]$, the theorem follows easily from Lemma \ref{lem:errorontriangles} and a combination of Lemma \ref{lem:uPiuonT} with Lemmas \ref{lem:upiuphi3} and \ref{lem:PiTPi}. \qed
\end{proof}

\begin{theorem}
	\label{teo:main}
	Let $K$ be an arbitrary convex quadrilateral satisfying the $MAC$ and the $clos$ $($eq. the $DAC$ and the $clos)$. If $l_1$ and $l_2$ are two neighboring sides of $K$ such that the parallelogram determined by them contains $K$ entirely, then, for any $p \in [1,\infty)$, we have
	\bq
	\label{eq:mainteo}
	|u-Qu|_{1,p,K} \le C \left[ 
	|l_1| \left\| \partial_{l_1} \nabla u \right\|_{0,p,K} + 
	|l_2| \left\| \partial_{l_2} \nabla u \right\|_{0,p,K} \right]
	\eq
	for some positive constant $C$ depending only on $p$ and those constants involved in the $MAC$ and the $clos$.
\end{theorem}

\begin{proof}
	Theorem \ref{teo:MAC+clos} ensures the existence of a reference element $K(a,b,\at,\bt)$ verifying conditions $[D1',D2]$ which is equivalent to $K$. In Section \ref{sub:red} we proved that if the anisotropic interpolation error estimate \eqref{eq:teoref} holds on $K(a,b,\at,\bt)$, then the anisotropic interpolation error estimate \eqref{eq:mainteo} holds on $K$. Thanks to Theorem \ref{teo:ref}, \eqref{eq:teoref} is valid and hence the proof is complete. \qed   
\end{proof}

\medskip
Finally, Counterexample 6.1 in \cite{AM1} shows that the interpolation error estimate
\bq
\label{eq:clas}
|u-Qu|_{1,p,K} \le Ch|u|_{2,p,K}
\eq
does not holds for $K=K(1,1,s,s)$, with $\frac{1}{2} <s< 1$, and $u(x)=x_1x_2$ when $p \ge 3$. It is clear that the anisotropic interpolation error estimate \eqref{eq:mainteo} implies \eqref{eq:clas}, so we conclude that \eqref{eq:mainteo} does not hold for the same election of $K$, $u$ and $p$. Since $K$ verifies the $clos$ but it does not obey the $MAC$ when $s$ tends to $1/2$, we conclude that the assumption $MAC$ in Theorem \ref{teo:main} can not to be relaxed. The question about whether the assumption $clos$ in Theorem \ref{teo:main} is also necessary is open.

\end{document}